\documentclass[12pt]{article}
\usepackage{amsfonts}
\usepackage{bbm}
\usepackage{mathrsfs}
\usepackage{amssymb}

\title{A Constructive Characterization of Solvable  \\
Polynomial Algebras
\thanks{Project supported by the National Natural Science
Foundation of China (10971044).}}

\author{Huishi Li\\
{\small Department of Applied Mathematics, College of Information Science and Technology}\\
{\small Hainan University,  Haikou 570228, China}}

\date{}

\begin{document}
\maketitle
\begin{center}
\begin{minipage}{120mm}
{\small {\bf Abstract.} For the solvable polynomial algebras
introduced and studied by Kandri-Rody and Weispfenning [J. Symbolic
Comput., 9(1990)],  a constructive characterization is given in
terms of Gr\"obner bases for ideals of free algebras, thereby
solvable polynomial algebras  are completely determinable in a
computational way.
 }\end{minipage}\end{center} {\parindent=0pt\vskip 6pt

\vskip 6pt

}

\def\QED{\hfill{$\Box$}}
\def \r{\rightarrow}
\def\T#1{\widetilde #1}
\def\NZ{\mathbb{N}}
\def\normalbaselines{\baselineskip 24pt\lineskip 4pt\lineskiplimit 4pt}
\def\mapdown#1{\llap{$\vcenter {\hbox {$\scriptstyle #1$}}$}
                                \Bigg\downarrow}
\def\mapdownr#1{\Bigg\downarrow\rlap{$\vcenter{\hbox
                                    {$\scriptstyle #1$}}$}}
\def\mapright#1#2{\smash{\mathop{\longrightarrow}\limits^{#1}_{#2}}}

\def\v5{\vskip .5truecm}

\def\OV#1{\overline {#1}}
\def\hang{\hangindent\parindent}
\def\textindent#1{\indent\llap{#1\enspace}\ignorespaces}
\def\item{\par\hang\textindent}


\def\LH{{\bf LH}}\def\LM{{\bf LM}}\def\LT{{\bf
LT}}\def\KX{K[\mathbbm{x}]} \def\KS{K\langle X\rangle}
\def\B{{\cal B}} \def\LC{{\bf LC}} \def\G{{\cal G}} \def\FRAC#1#2{\displaystyle{\frac{#1}{#2}}}
\def\SUM^#1_#2{\displaystyle{\sum^{#1}_{#2}}} \def\O{{\cal O}}  \def\J{{\bf J}}\def\BE{\B (\mathbbm{e})}
\def\PRCVE{\prec_{\varepsilon\hbox{-}gr}}\def\BV{\B (\varepsilon )}\def\PRCEGR{\prec_{\mathbbm{e}\hbox{\rm -}gr}}
\def\KM{K[\mathbbm{x}]}\def\H{\cal H}

\def\KS{K\langle X\rangle}
\def\LR{\langle X\rangle}
\def\HL{{\rm LH}}
\vskip 1truecm

\section*{1. Introduction}

In the late 1980s, the celebrated  Gr\"obner basis theory developed
by Bruno Buchberger [Bu1, 2] for commutative polynomial ideals  was
successfully generalized to one-sided ideals in enveloping algebras
of Lie algebras by Apel and Lassner [AL], to one-sided ideals in
Weyl algebras (including  algebras of partial differential operators
with polynomial coefficients over a field of characteristic 0) by
Galligo [Gal], and more generally, to one-sided and two-sided ideals
in  solvable polynomial algebras (or algebras of solvable type) by
Kandri-Rody and Weispfenning [K-RW]. The class of solvable
polynomial algebras includes  not only the commutative polynomial
algebras, the Weyl algebras, the enveloping algebras of finite
dimensional Lie algebras, and a large number of iterated Ore
extensions, but also numerous other significant noncommutative
algebras (cf. [Li1], [BGV], [Li4]).
\par

Let $K$ be a field and $R=K[X_1,\ldots ,X_n]$ the commutative
polynomial $K$-algebra in $n$ variables. Originally, a
noncommutative solvable polynomial algebra $R'$ was defined in
[K-RW] by first fixing a monomial ordering $\prec$ on the standard
$K$-basis $\mathscr{B} =\{ X_1^{\alpha_1}\cdots
X_n^{\alpha_n}~|~\alpha_i\in\NZ\}$ of $R$, and then introducing a
new multiplication $*$ on $R$, such that certain axioms ([K-RW],
AXIOMS 1.2) are satisfied. In the formal language of associative
$K$-algebras, a solvable polynomial algebra can actually be defined
as a finitely generated associative $K$-algebra $A=K[a_1,\ldots
,a_n]$, that has the PBW $K$-basis $\B =\{ a_1^{\alpha_1}\cdots
a_n^{\alpha_n}~|~\alpha_i\in\NZ\}$ and a (two-sided) monomial
ordering $\prec$ on $\B$ such that for $1\le i<j\le n$,
$a_ja_i=\lambda_{ji}a_ia_j+f_{ji}$ and $\LM (f_{ji})\prec a_ia_j$,
where $\lambda_{ji}\in K-\{ 0\}$ and  $f_{ji}\in K$-span$\B$ ([LW],
Definition 2.1). Full details on this definition will be recalled in
the next section. \par

Let $\KS =K\langle X_1,\ldots ,X_n\rangle$ be the noncommutative
free $K$-algebra on $X=\{ X_1,\ldots ,X_n\}$, and  $\mathbb{B}=\{
1,~X_{i_1}X_{i_2}\cdots X_{i_s}~|~X_{i_j}\in X,~s\ge 1\}$ the
standard $K$-basis of $\KS$. For convenience, we use capital letters
$U,V,W,S,\ldots$ to denote elements (monomials) of $\mathbb{B}$.
Recall that a monomial ordering $\prec '$ on $\mathbb{B}$ is a
well-ordering such that for $W,U,V\in\mathbb{B}$, $U\prec 'V$
implies $WU\prec 'WV$, $UW\prec 'VW$; and moreover, if
$W,U,V,S\in\mathbb{B}$ with $W\ne V$, then $W=UVS$ implies $V\prec
'W$ (thereby $1\prec 'W$ for all $1\ne W\in\mathbb{B}$). As before,
we let  $R=K[X_1,\ldots X_n]$ denote the commutative polynomial
algebra in variables $X_1,\ldots X_n$, and write $\mathscr{B} =\{
X_1^{\alpha_1}\cdots X_n^{\alpha_n}~|~\alpha_i\in\NZ\}$ for the
standard $K$-basis of $R$. In [K-RW] it was shown how to  construct
a solvable polynomial algebra by means of a fixed monomial ordering
$\prec$ on  $\mathscr{B}$ and a {\it commutation system} $Q=\{
q_{ij}~|~1\le i<j\le n\}\subset \KS$, where each
$q_{ij}=X_jX_i-c_{ij}X_iX_j-p_{ij}$ with $c_{ij}\in K-\{ 0\}$ and
$p_{ij}\in R$ such that $p_{ij}\prec X_iX_j$, that is,
{\parindent=.5truecm \par

\item{$\bullet$} ([K-RW], Theorem 1.7) if $I(Q)$
denotes the (two-sided) ideal of $\KS$ generated by $Q$, then the
quotient algebra $\KS /I(Q)$ is isomorphic to a solvable polynomial
algebra $R'$ defined with respect to the given $\prec$ on
$\mathscr{B}$ and the new multiplication $*$ on $R$ satisfying
([K-RW], AXIOMS 1.2) if and only if $I(Q)$ satisfies the condition
$${\rm (H)}\hskip 3truecm I(Q)~\hbox{contains no nonzero commutative
polynomials} .$$\par {\parindent=0pt

With a fixed monomial ordering $\prec$ on $\mathscr{B}$ and a given
commutation system $Q$ as above, the relation between the condition
(H) and Gr\"obner bases of $I(Q)$ in $\KS$ was also explored in
[K-RW], that is, {\parindent=.5truecm\par

\item{$\bullet$} ([K-RW], Theorem 1.11 due to Mora) assuming that there is a positive
monomial ordering $\prec '$ on
$\mathbb{B}$ which extends the given monomial ordering $\prec$ on
$\mathscr{B}$, such that
$$X_1^{\alpha_1}X_2^{\alpha_2}\cdots
X_j^{\alpha_j}\prec 'X_jX_i~\hbox{for all monomials}~
X_1^{\alpha_1}X_2^{\alpha_2}\cdots X_j^{\alpha_j},~1\le i<j\le n,$$
then $I(Q)$ satisfies the condition (H) if and only if $Q$ is a
Gr\"obner basis for $I(Q)$ with respect to the monomial ordering
$\prec '$ on $\mathbb{B}$.\par {\parindent=0pt

As an example, in the case that $\prec$ is the pure lexicographic
ordering on $\B$, the existence of such a positive monomial ordering
$\prec '$ on $\mathbb{B}$ was given by Mora (see the ordering $\prec
'$ described in [K-RW] after Corollary 1.12).}\par

In the case that $Q=\{ X_jX_i-X_iX_j-p_{ij}(X_1,\ldots
,X_{j-1})~|~1\le i<j\le n\}\subset\KS$ but with each
$p_{ij}(X_1,\ldots ,X_{j-1})$ being an element of the free
$K$-algebra $K\langle X_1,\ldots ,X_{j-1}\rangle$ on $\{ X_1,\ldots
,X_{j-1}\}$, it was proved that {\parindent=.5truecm \par

\item{$\bullet$} ([K-RW], Theorem 1.13) the quotient algebra $\KS /I(Q)$ is a solvable polynomial
algebra of strictly lexicographical type if and only if $I(Q)$
satisfies the condition (H).\v5

In [Li1, 4], some results on the construction of solvable polynomial
algebras by means of Gr\"obner bases in $\KS$ were also given
([Li1], CH.III, Proposition 2.2,  Proposition 2.3; [Li4], Ch.4,
Proposition 4.2). However, so far there seems no a complete
constructive characterization of solvable polynomial algebras in
terms of Gr\"obner bases for ideals of free algebras,  by which one
may effectively determine all such algebras. In this note, we solve
this problem in Section 2 by Theorem 2.5. Furthermore, we remark in
Section 3 that the main results obtained for quadric solvable
polynomial algebras in [Li1] indeed hold true for arbitrary solvable
polynomial algebras; and  we remark that the algebras satisfying the
equivalent conditions of Proposition 2.4, especially the algebras
described in Theorem 2.5(ii)(a),  also provide us with an
interesting class of algebras in the computational noncommutative
algebra. \v5

Throughout  this paper, $K$ denotes a field, $K^*=K-\{ 0\}$;
$\mathbb{N}$ denotes the set of all nonnegative integers.
Moreover, the Gr\"obner basis theory for ideals of free algebras is
referred to [Mor].

\section*{2. The Main Result}
We first briefly recall from ([K-RW], [LW], [Li1]) some basics
concerning solvable polynomial algebras. Let $A=K[a_1,\ldots ,a_n]$
be a finitely generated $K$-algebra with the minimal set of
generators $\{ a_1,\ldots ,a_n\}$. If, for some permutation $\tau
=i_1i_2\cdots i_n$ of $1,2,\ldots ,n$, the set $\B =\{
a^{\alpha}=a_{i_1}^{\alpha_1}\cdots a_{i_n}^{\alpha_n}~|~\alpha
=(\alpha_1,\ldots ,\alpha_n)\in\NZ^n\} ,$ forms a $K$-basis of $A$,
then $\B$ is referred to as a {\it PBW $K$-basis} of $A$. It is
clear that if $A$ has a PBW $K$-basis, then we can always assume
that $i_1=1,\ldots ,i_n=n$. Thus, we make the following convention
once for all.{\parindent=0pt\v5

{\bf Convention} From now on in this paper, if we say that the
algebra $A$ has the PBW $K$-basis $\B$, then $\B$ is meant the one
$$\B =\{ a^{\alpha}=a_{1}^{\alpha_1}\cdots
a_{n}^{\alpha_n}~|~\alpha =(\alpha_1,\ldots ,\alpha_n)\in\NZ^n\} .$$
Moreover, adopting the commonly used terminology in computational
algebra, elements of $\B$ are referred to as {\it monomials} of
$A$.} \v5

Suppose that $A$ has the PBW $K$-basis $\B$ as presented above and
that  $\prec$ is a total ordering on $\B$. Then every nonzero
element $f\in A$ has a unique expression
$$f=\lambda_1a^{\alpha (1)}+\lambda_2a^{\alpha (2)}+\cdots +\lambda_ma^{\alpha (m)},~
\lambda_j\in K^*,~a^{\alpha
(j)}=a_1^{\alpha_{1j}}a_2^{\alpha_{2j}}\cdots a_n^{\alpha_{nj}}\in\B
,~1\le j\le m,$$  in which the basis elements satisfy  $a^{\alpha
(1)}\prec a^{\alpha (2)}\prec\cdots \prec a^{\alpha (m)}$.  It
follows that we may define
$$\begin{array}{ll} \LM (f)=a^{\alpha (m)},&\hbox{the leading monomial of}~f;\\
\LC (f)=\lambda_m,&\hbox{the leading coefficient of}~f;\\
\LT (f)=\lambda_ma^{\alpha (m)},&\hbox{the leading term
of}~f.\end{array}$${\parindent=0pt\v5

{\bf 2.1. Definition}  Suppose that the $K$-algebra $A=K[a_1,\ldots
,a_n]$ has the PBW $K$-basis $\B$. If $\prec$ is a total ordering on
$\B$ that satisfies the following three conditions:\par

(1) $\prec$ is a well-ordering;\par

(2) For $a^{\gamma},a^{\alpha},a^{\beta}\in\B$, if $a^{\alpha}\prec
a^{\beta}$ and $\LM (a^{\gamma}a^{\alpha})$, $\LM
(a^{\gamma}a^{\beta})\not\in K$, then $\LM
(a^{\gamma}a^{\alpha})\prec\LM (a^{\gamma}a^{\beta})$;\par

(3) For $a^{\gamma},a^{\alpha},a^{\beta}\in\B$, if $a^{\beta}\ne
a^{\gamma}$, and $a^{\gamma}=\LM (a^{\alpha}a^{\beta})$, then
$a^{\beta}\prec a^{\gamma}$ (thereby $1\prec a^{\gamma}$ for all
$a^{\gamma}\ne 1$),\par

then we call $\prec$ a {\it left monomial ordering} on $\B$. }\par

Similarly, a {\it right monomial ordering} on $\B$ may be
defiend.\par

If $\prec$ is both a left monomial ordering and a right monomial
ordering on $\B$, then we call $\prec$ a {\it two-sided monomial
ordering} (or simply a {\it monomial ordering}) on $\B$.
{\parindent=0pt\v5

{\bf Remark} Examples of left (right) monomial orderings on PBW
$K$-bases, which are not right (left) monomial orderings, are given
in [Li3].\v5

{\bf 2.2. Definition} Suppose that the $K$-algebra $A=K[a_1,\ldots
,a_n]$ has the PBW $K$-basis $\B$ and that $\prec$ is a (two-sided)
monomial ordering on $\B$. If for all
$a^{\alpha}=a_1^{\alpha_1}\cdots a_n^{\alpha_n}$,
$a^{\beta}=a_1^{\beta_1}\cdots a^{\beta_n}_n\in\B$, the following
holds:
$$\begin{array}{rcl} a^{\alpha}a^{\beta}&=&\lambda_{\alpha ,\beta}a^{\alpha +\beta}+f_{\alpha ,\beta},\\
&{~}&\hbox{where}~\lambda_{\alpha ,\beta}\in K^*,~a^{\alpha
+\beta}=a_1^{\alpha_1+\beta_1}\cdots a_n^{\alpha_n+\beta_n},~\hbox{and either}~f_{\alpha ,\beta}=0~\hbox{or}\\
&{~}&f_{\alpha ,\beta}\in K\hbox{-span}\B~\hbox{with}~\LM (f_{\alpha
,\beta})\prec a^{\alpha +\beta},\end{array}$$ then $A$ is said to be
a {\it solvable polynomial algebra}.} \v5

The results of the next proposition are summarized from ([K-RW],
Sections 2 -- 5).{\parindent=0pt\v5

{\bf 2.3. Proposition}  Let $A=K[a_1,\ldots ,a_n]$ be a solvable
polynomial algebra with the (two-sided) monomial ordering $\prec$ on
the PBW $K$-basis $\B$ of $A$. The following statements hold.\par

(i) $A$ is a (left and right) Noetherian domain.\par

(ii) Every left ideal $I$ of $A$ has a finite  left Gr\"obner basis
$\G=\{ g_1,\ldots ,g_t\}$ in the sense that}{\parindent=.5truecm\par

\item{$\bullet$} if $0\ne f\in I$, then there is some $g_i\in\G$ such that $\LM (g_i)|\LM (f)$,
i.e., there is some $a^{\gamma}\in\B$ such that $\LM (f)=\LM
(a^{\gamma}\LM (g_i))$, or equivalently, with $\gamma
(i_j)=(\gamma_{i_{1j}},\gamma_{i_{2j}},\ldots
,\gamma_{i_{nj}})\in\NZ^n$, $f$ has a left Gr\"obner representation:
$$\begin{array}{rcl} f&=&\sum_{i,j}\lambda_{ij}a^{\gamma (i_j)}g_j,~\hbox{where}~\lambda_{ij}\in K^*,
~a^{\gamma (i_j)}\in\B ,~
g_j\in \G,\\
&{~}&\hbox{satisfying}~\LM (a^{\gamma (i_j)}g_j)\preceq\LM
(f)~\hbox{for all}~(i,j).\end{array}$$}{\parindent=0pt\par

(iii) The Buchberger's Algorithm, that computes a finite Gr\"obner
basis for a finitely generated commutative polynomial ideal, has a
complete noncommutative version that computes a finite left
Gr\"obner basis for a finitely generated left ideal
$I=\sum_{i=1}^mAf_i$ of $A$.\par

(iv) Similar results of (ii) and (iii) hold for right ideals and
two-sided ideals of $A$.\par\QED}\v5

It is clear that a solvable polynomial algebra $A$ depends on two
independent factors, namely a PBW $K$-basis $\B$ and an appropriate 
monomial ordering $\prec$ on $\B$. Let $\KS =K\langle X_1,\ldots 
,X_n\rangle$ be the free $K$-algebra on $X=\{ X_1,\ldots ,X_n\}$ and 
$\mathbb{B}=\{ 1,~X_{i_1}\cdots X_{i_s}~|~X_{i_j}\in X,~s\ge 1\}$ 
the standard $K$-basis of $\KS$. Let $I$ be an ideal of $\KS$. 
Concerning the relation between Gr\"obner bases of $I$ and the 
existence of a PBW $K$-basis for  the quotient algebra $A=\KS /I$, 
we recall  the following{\parindent=0pt\v5

{\bf 2.4. Proposition} ([Li4], Ch 4, Theorem 3.1) Let $A=\KS /I$ be
as above. Suppose that $I$ contains a subset of $\frac{n(n-1)}{2}$
elements
$$G=\{ g_{ji}=X_jX_i-F_{ji}~|~F_{ji}\in\KS ,~1\le i<j\le n\}$$
such that with respect to some monomial ordering $\prec_{_X}$ on
$\mathbb{B}$, $\LM (g_{ji})=X_jX_i$ holds for all the $g_{ji}$. The
following two statements are equivalent:\par

(i) $A$ has the  PBW $K$-basis $\mathscr{B}=\{ \OV X_1^{\alpha_1}\OV
X_2^{\alpha_2}\cdots \OV X_n^{\alpha_n}~|~\alpha_j\in\NZ\}$ where
each $\OV X_i$ denotes the coset of $I$ represented by $X_i$ in $\KS
/I$.\par

(ii) Any  subset $\G$ of $I$ containing $G$ is a  Gr\"obner basis
for $I$ with respect to $\prec_{_X}$.\QED \v5

{\bf Remark} Obviously, Proposition 2.4 holds true if we use any
permutation $\{X_{k_1},\ldots ,X_{k_n}\}$ of $\{ X_1,\ldots ,X_n\}$.
So, in what follows we conventionally use only $\{ X_1,\ldots
,X_n\}$.}\v5

We note that that if $G=\{ g_{ji}=X_jX_i-F_{ji}~|~F_{ji}\in\KS
,~1\le i<j\le n\}$ is a Gr\"obner basis of the ideal $I$ such that
$\LM (g_{ji})=X_jX_i$ for all the $g_{ji}$, then the {\it reduced
Gr\"obner basis} of $I$ is of the form
$$\G =\left\{\left.
g_{ji}=X_jX_i-\sum_q\mu^{ji}_qX_1^{\alpha_{1q}}X_2^{\alpha_{2q}}\cdots
X_n^{\alpha_{nq}}~\right |~ \LM (g_{ji})=X_jX_i,~ 1\le i<j\le n
\right\}
$$  where $\mu^{ji}_q\in K$ and $(\alpha_{1q},\alpha_{2q},\ldots ,\alpha_{nq})\in\NZ^n$.
Bearing in mind Definition 2.2 and combining this fact, we have the
following characterization of solvable polynomial algebras in terms
of Gr\"obner bases for ideals of free algebras. {\parindent=0pt\v5

{\bf 2.5. Theorem}  Let $A=K[a_1,\ldots ,a_n]$ be a finitely
generated algebra over the field $K$, and let $\KS =K\langle
X_1,\ldots ,X_n\rangle$ be the free $K$-algebras with the standard
$K$-basis $\mathbb{B}=\{ 1,~X_{i_1}\cdots X_{i_s}~|~X_{i_j}\in
X,~s\ge 1\}$. With notation as before, the following two statements
are equivalent:\par

(i) $A$ is a solvable polynomial algebra in the sense of Definition
2.2.\par

(ii) $A\cong \OV A=\KS /I$ via the $K$-algebra epimorphism $\pi_1$:
$\KS \r A$ with $\pi_1(X_i)=a_i$, $1\le i\le n$, $I=$ Ker$\pi_1$,
satisfying  {\parindent=1.3truecm

\item{(a)} with respect to some monomial ordering $\prec_{_X}$ on $\mathbb{B}$, the ideal $I$ has a
finite Gr\"obner basis $G$ and the reduced Gr\"obner basis of $I$ is
of the form
$$\G =\left\{\begin{array}{rcl} g_{ji}&=&X_jX_i-\lambda_{ji}X_iX_j-F_{ji}\\
&{~}&\hbox{with}~F_{ji}=\sum_q\mu^{ji}_qX_1^{\alpha_{1q}}X_2^{\alpha_{2q}}\cdots 
X_n^{\alpha_{nq}}\end{array}~\left |~\begin{array}{l} \LM 
(g_{ji})=X_jX_i,\\ 1\le i<j\le n\end{array}\right. \right\} $$ where 
$\lambda_{ji}\in K^*$, $\mu^{ji}_q\in K$, and 
$(\alpha_{1q},\alpha_{2q},\ldots ,\alpha_{nq})\in\NZ^n$, thereby   
$\mathscr{B}=\{ \OV X_1^{\alpha_1}\OV X_2^{\alpha_2}\cdots \OV 
X_n^{\alpha_n}~|~$ $\alpha_j\in\NZ\}$ forms  a PBW $K$-basis for 
$\OV A$, where each $\OV X_i$ denotes the coset of $I$ represented 
by $X_i$ in $\OV A$; and

\item{(b)} there is a (two-sided) monomial ordering
$\prec$ on $\mathscr{B}$ such that $\LM (\OV{F}_{ji})\prec \OV
X_i\OV X_j$ whenever $\OV F_{ji}\ne 0$, where $\OV
F_{ji}=\sum_q\mu^{ji}_q\OV X_1^{\alpha_{1i}}\OV
X_2^{\alpha_{2i}}\cdots \OV X_n^{\alpha_{ni}}$, $1\le i<j\le n$.
\vskip 6pt}

{\bf Proof} (i) $\Rightarrow$ (ii) Let $\B=\{
a^{\alpha}=a_1^{\alpha_1}\cdots a_n^{\alpha_n}~|~\alpha
=(\alpha_1,\ldots ,\alpha_n)\in\NZ^n\}$ be the PBW $K$-basis of the
solvable polynomial algebra $A$ and $\prec$ a (two-sided) monomial
ordering on $\B$. By Definition 2.2, the generators of $A$ satisfy
the relations:
$$a_ja_i=\lambda_{ji}a_ia_j+f_{ji},\quad 1\le i<j\le n,\eqno{(*)}$$
where $\lambda_{ji}\in K^*$ and  $f_{ji}=\sum_q\mu^{ji}_qa^{\alpha
(q)}\in K$-span$\B$ with $\LM (f_{ji})\prec a_ia_j$. Consider in the
free $K$-algebra $\KS =K\langle X_1,\ldots ,X_n\rangle$ the subset
$$\G =\{ g_{ji}=X_jX_i-\lambda_{ji}X_iX_j-F_{ji}~|~1\le i<j\le n\} ,$$
where if
$f_{ji}=\sum_q\mu^{ji}_qa_1^{\alpha_{1q}}a_2^{\alpha_{2q}}\cdots
a_n^{\alpha_{nq}}$ then
$F_{ji}=\sum_q\mu^{ji}_qX_1^{\alpha_{1q}}X_2^{\alpha_{2q}}\cdots
X_n^{\alpha_{nq}}$ for $1\le i<j\le n$. We write
$J=\langle\G\rangle$ for the ideal of $\KS$ generated by $\G$ and
put $\OV A=\KS /J$. Let $\pi_1$: $\KS\r A$ be the $K$-algebra
epimorphism with $\pi_1(X_i)=a_i$, $1\le i\le n$, and let $\pi_2$:
$\KS\r\OV A$ be the canonical algebra epimorphism. It follows from
the fundamental theorem of homomorphism that there is an algebra
epimorphism $\varphi$: $\OV A\r A$ defined by $\varphi (\OV
X_i)=a_i$, $1\le i\le n$, such that the following diagram of algebra
homomorphisms is commutative:
$$\begin{array}{cccc}
\KS&\mapright{\pi_2}{}&\OV A&\\
\mapdown{\pi_1}&\swarrow\scriptstyle{\varphi}&&\varphi\circ\pi_2=\pi_1\\
A&&&\end{array}$$  On the other hand, by the definition of each
$g_{ji}$ we see that every element $\OV H\in\OV A$ may be written as
$\OV H=\sum_j\mu_j\OV X_1^{\beta_{1j}}\OV X_2^{\beta_{2j}}\cdots \OV
X_n^{\beta_{nj}}$ with $\mu_j\in K$ and $(\beta_{1j},\ldots
,\beta_{nj})\in\NZ^n$, where each $\OV X_i$ is the coset of $J$
represented by $X_i$ in $\OV A$. Noticing the relations presented in
$(*)$,  it is straightforward to check that the correspondence
$$\begin{array}{cccc} \psi :&A&\mapright{}{}&\OV A\\
&\displaystyle{\sum_i}\lambda_ia_1^{\alpha_{1i}}\cdots
a_n^{\alpha_{ni}}&\mapsto&\displaystyle{\sum_i}\lambda_i\OV
X_1^{\alpha_{1i}}\cdots \OV X_n^{\alpha_{ni}}\end{array}$$ is an
algebra homomorphism such that $\varphi\circ\psi =1_A$ and
$\psi\circ\varphi =1_{\OV A}$, where $1_A$ and $1_{\OV A}$ denote
the multiplicative identities of $A$ and $\OV A$ respectively. This
shows that $A\cong\OV A$, thereby Ker$\pi_1=I=J$; moreover,
$\mathscr{B}=\{ \OV X_1^{\alpha_1}\OV X_2^{\alpha_2}\cdots \OV
X_n^{\alpha_n}~|~\alpha_j\in\NZ\}$ forms  a PBW $K$-basis for $\OV
A$, and $\prec$ is a (two-sided) monomial ordering on
$\mathscr{B}$.}
\par
We next show that $\G$ forms the reduced Gr\"obner basis for $I$ as
described in (a). To this end, we first show that the monomial
ordering $\prec$ on $\B$ induces a monomial ordering $\prec_{_X}$ on
the standard $K$-basis $\mathbb{B}$ of $\KS$. For convenience, we
use capital letters $U,V,W,S,\ldots$ to denote elements (monomials)
in $\mathbb{B}$. We also fix a graded lexicographic ordering
$\prec_{grlex}$ on $\mathbb{B}$ (with respect to a fixed positively
weighted gradation of $\KS$) such that
$$X_1\prec_{grlex}X_2\prec_{grlex}\cdots\prec_{grlex}X_n.$$
Then, for $U,V\in\mathbb{B}$ we define
$$U\prec_{_X}V~\hbox{if}~\left\{\begin{array}{l} \LM (\pi_1 (U))\prec\LM (\pi_1(V)),\\
\hbox{or}\\
\LM (\pi_1 (U))=\LM (\pi_1
(V))~\hbox{and}~U\prec_{grlex}V.\end{array}\right.$$ Since  $A$ is a
domain (Proposition 2.3(i)) and $\pi_1$ is an algebra homomorphism
with $\pi_1(X_i)=a_i$ for $1\le i\le n$, it follows that $\LM
(\pi_1(W))\ne 0$ for all $W\in\mathbb{B}$. We also note from
Definition 2.2 that if $f,g\in A$ are nonzero elements, then $\LM
(fg)=\LM (\LM (f)\LM (g))$. Thus, if $U,V,W\in\mathbb{B}$ and
$U\prec_{_X}V$ subject to $\LM (\pi_1(U))\prec\LM (\pi_1(V))$, then
$$\begin{array}{cc}\begin{array}{rcl} \LM
(\pi_1(WU))&=&\LM (\LM (\pi_1(W))\LM (\pi_1(U)))\\
&\prec&\LM (\LM (\pi_1(W))\LM (\pi_1(V)))\\
&=&\LM (\pi_1(WV))\end{array}&\end{array}$$ implies
$WU\prec_{_X}WV$; if $U\prec_{_X}V$ subject to $\LM (\pi_1 (U))=\LM
(\pi_1 (V))~\hbox{and}~U\prec_{grlex}V$, then
$$\begin{array}{rcl} \LM (\pi_1(WU))&=&\LM (\LM (\pi_1(W))\LM (\pi_1(U)))\\
&=&\LM (\LM (\pi_1(W))\LM (\pi_1(V)))\\
&=&\LM (\pi_1(WV))\end{array}$$ \\
and $WU\prec_{grlex}WV$ implies $WU\prec_{_X}WV$.  Similarly, if 
$U\prec_{_X}V$ then $US\prec_{_X}VS$ for all $S\in\mathbb{B}$. 
Moreover, if $W,U,V,S\in\mathbb{B}$, $W\ne V$, such that $W=UVS$, 
then $\LM (\pi_1(W))=\LM (\pi_1(UVS))$ and clearly 
$V\prec_{grlex}W$, thereby $V\prec_{_X}W$. Since  $\prec$ is a 
well-ordering on $\B$ and $\prec_{grlex}$ is a well-ordering on 
$\mathbb{B}$, the above argument shows that $\prec_{_X}$ is a 
(two-sided) monomial ordering on $\mathbb{B}$. With this monomial 
ordering $\prec_{_X}$ in hand,  by the definition of $F_{ji}$ we see 
that $\LM (F_{ji})\prec_{_X}X_iX_j$. Furthermore, since $\LM 
(\pi_1(X_jX_i))=a_ia_j=\LM (\pi_1(X_iX_j))$ and $X_iX_j\prec_{grlex} 
X_jX_i$, we see that $X_iX_j\prec_{_X}X_jX_i$. It follows that $\LM 
(g_{ji})=X_jX_i$ for $1\le i<j\le n$. Now, by Proposition 2.4 we 
conclude that $\G$ forms a Gr\"obner basis for $I$ with respect to 
$\prec_{_X}$. Finally, by the definition of $\G$, it is clear that 
$\G$ is the reduced Gr\"obner basis of $I$ and $\G$ meets the 
requirement of Proposition 2.4,  hence $\mathscr{B}=\{ \OV 
X_1^{\alpha_1}\OV X_2^{\alpha_2}\cdots \OV 
X_n^{\alpha_n}~|~\alpha_j\in\NZ\}$ forms a PBW $K$-basis for $\OV 
A$, as desired.\par

(ii) $\Rightarrow$ (i)  Note that (a) $+$ (b) tells us that  the 
generators of $\OV A$ satisfy the relations $\OV X_j\OV 
X_i=\lambda_{ji}\OV X_i\OV X_j+\OV F_{ji}$, $1\le i<j\le n$, and 
that if $\OV F_{ji}\ne 0$ then $\LM (\OV F_{ji})\prec \OV X_i\OV 
X_j$ with respect to the given monomial ordering $\prec$ on 
$\mathscr{B}$. It follows that  $\OV A$ and hence $A$ is a solvable 
polynomial algebra in the sense of Definition 2.2.\QED 
{\parindent=0pt\v5

{\bf Remark} The monomial ordering $\prec_{_X}$ we defined in the
proof of Theorem 2.5 is a modification of the {\it lexicographic
extension} defined in [EPS]. But our definition of $\prec_{_X}$
involves a graded monomial ordering $\prec_{grlex}$ on the standard
$K$-basis $\mathbb{B}$ of the free $K$-algebra $\KS =K\langle
X_1,\ldots ,X_n\rangle$. The reason is that {\it the monomial
ordering $\prec_{_X}$ on $\mathbb{B}$ must be compatible with the
usual rule of division}, namely, $W,U,V,S\in\mathbb{B}$, $W\ne V$,
and $W=UVS$ implies $V\prec_{_X}W$. While it is clear that if we use
any lexicographic ordering $\prec_{lex}$ in the definition of
$\prec_{_X}$, then this rule will not work in general.}\v5

Note that our discussion on solvable polynomial algebras made use of
Definition 2.1 for a monomial ordering and Definition 2.2 for a
solvable polynomial algebra. In light of ([K-RW], AXIOMS 1.2), the
next proposition may make the practical use of Theorem 2.5 much
easier and flexible in determining solvable polynomial algebras.
{\parindent=0pt\v5

{\bf 2.6. Proposition} Let $A=K[a_1,\ldots ,a_n]$  be a finitely
$K$-algebra with the PBW $K$-basis $\B=\{
a^{\alpha}=a_1^{\alpha_1}\cdots a_n^{\alpha_n}~|~\alpha
=(\alpha_1,\ldots ,\alpha_n)\in\NZ^n\}$ and the generators of $A$
satisfy the relations:
$$a_ja_i=\lambda_{ji}a_ia_j+f_{ji},\quad 1\le i<j\le n,$$
where $\lambda_{ji}\in K^*$ and $f_{ji}\in K$-span$\B$. Then the
following two statements are equivalent.\par

(i) There is a (two-sided) monomial ordering $\prec$ on $\B$ in the
sense of Definition 2.1, such that
$$\LM (f_{ji})\prec a_ia_j,\quad 1\le i<j\le n.$$\par

(ii) There is a monomial ordering $\prec$ on the additive monoid
$\NZ^n$, i.e., $\prec$ is a well-ordering and $\alpha\prec\beta$
implies $\alpha +\gamma\prec\beta +\gamma$ for all $\alpha ,\beta
,\gamma\in\NZ^n$, such that if $\LM (f_{ji})=a^{\alpha}$ and we
write $a^{\beta}=a_ia_j$, then $\alpha\prec\beta$, $1\le i<j\le
n$.\par\QED}

We end this section by an example illustrating Theorem 2.5 and
Proposition 2.6.{\parindent=0pt\v5

{\bf Example 1.}  Considering the $\NZ$-graded structure of the free
$K$-algebra $\KS =K\langle X_1,X_2, X_3\rangle$ by assigning $X_1$
the degree 2, $X_2$ the degree 1 and $X_3$ the degree 4, let $I$ be
the ideal of $\KS$ generated by the elements
$$\begin{array}{l} g_1=X_1X_2- X_2X_1,\\
g_2=X_3X_1-\lambda X_1X_3-\mu X_3X_2^2-f(X_2),\\
g_3=X_3X_2- X_2X_3,\end{array}$$ where $\lambda\in K^*$, $\mu\in K$,
$f(X_2)$ is a polynomial in $X_2$ which has degree $\le 6$, or
$f(X_2)=0$. The following properties hold.   \par

(1) If we use the graded lexicographic ordering
$X_2\prec_{grlex}X_1\prec_{grlex}X_3$ on $\KS$, then the three
generators have the leading monomial $\LM (g_1)=X_1X_2$, $\LM
(g_2)=X_3X_1$, and $\LM (g_3)=X_3X_2$. It is straightforward to
verify that $\G =\{ g_1,g_2,g_3\}$ forms a Gr\"obner basis for $I$.
\par

(2) With respect to the fixed $\prec_{grlex}$ in (1), the reduced
Gr\"obner basis $\G '$ of $I$ consists of
$$\begin{array}{l} g_1=X_1X_2- X_2X_1,\\
g_2=X_3X_1- \lambda X_1X_3-\mu X_2^2X_3-f(X_2),\\
g_3=X_3X_2- X_2X_3,\end{array}$$\par

(3) Writing $A=K[a_1,a_2, a_3]$ for the quotient algebra $\KS /I$,
where $a_1$, $a_2$ and $a_3$ denote the cosets $X_1+I$, $X_2+I$ and
$X_3+I$ in $\KS /I$ respectively, it follows that $A$ has the PBW
basis $\B =\{
a^{\alpha}=a_2^{\alpha_2}a_1^{\alpha_1}a_3^{\alpha_3}~|~\alpha
=(\alpha_2,\alpha_1,\alpha_3)\in\NZ^3\}$. Noticing that $
a_2a_1=a_1a_2$, it is clear that $\B ' =\{
a^{\alpha}=a_1^{\alpha_1}a_2^{\alpha_2}a_3^{\alpha_3}~|~\alpha
=(\alpha_1,\alpha_2,\alpha_3)\in\NZ^3\}$ is also a PBW basis for
$A$. Since $a_3a_1=\lambda a_1a_3+\mu a_2^2a_3+f(a_2)$, where
$f(a_2)\in K$-span$\{ 1,a_2,a_2^2,\ldots ,a_2^6\}$, we see that $A$
has the monomial ordering $\prec_{lex}$ on $\B '$, which is given by
the  lexicographic ordering $\prec_{lex}$ on $\NZ^3$ such that
$a_3\prec_{lex}a_2\prec_{lex}a_1$ and $\LM (\mu
a_2^2a_3+f(a_2))\prec_{lex}a_1a_3$, thereby $A$ is turned into a
solvable polynomial algebra with respect to $\prec_{lex}$.}\par

Moreover, one easily checks that if  $a_1$ is assigned the degree 2,
$a_2$ is assigned the degree 1 and $a_3$ is assigned the degree 4,
then, $A$ has another monomial ordering $\prec_{grlex}$ on $\B '$,
which is given by the graded lexicographic ordering $\prec_{grlex}$
on $\NZ^3$ such that  $a_3\prec_{grlex}a_2\prec_{grlex}a_1$ and $\LM
(\mu a_2^2a_3+f(a_2))\prec_{grlex}a_1a_3$, thereby $A$ is turned
into a solvable polynomial algebra with respect to $\prec_{grlex}$.
\v5

\section*{3. Further Remarks and Questions}

First recall from [Li1] that a {\it quadric solvable polynomial
algebra} is a solvable polynomial algebra $A=K[a_1,\ldots ,a_n]$
with the PBW $K$-basis $\B$ and a {\it graded monomial ordering}
$\prec_{gr}$ subject to the convention that {\it each $a_i$ is
assigned the degree 1}, such that
$$a_ja_i=\lambda_{ji}a_ia_j+\sum_{k\le \ell}\lambda_{ji}^{k\ell}a_ka_{\ell}+\sum_h\lambda_ha_h+c_{ji},
\quad 1\le i<j\le n,$$ where $\lambda_{ji}\in K^*$,
$\lambda_{ji}^{k\ell},\lambda_h,c_{ji}\in K$; in the case where
$\sum_{k\le
\ell}\lambda_{ji}^{k\ell}a_ka_{\ell}+\sum_h\lambda_ha_h+c_{ji}=0$
for $1\le i<j\le n$, a quadric solvable polynomial algebra $A$ is
called a {\it homogeneous solvable polynomial algebra}. For a
quadric solvable polynomial algebra $A$, by introducing  the
$\prec_{gr}$-filtration ${\cal F}A$ of $A$ with respect to
$\prec_{gr}$, and passing to the associated graded algebra $G^{\cal
F}(A)$ which is a homogeneous solvable polynomial algebra, ([Li1],
CH.V, CH.VI) shows how to calculate the Gelfand-Kirillov dimension
GK.dim$A/L$ and the multiplicity $e(A/L)$ of a cyclic $A$-module
$A/L$, as well as how to calculate GK.dim$(A/L\otimes_KA/J)$ and
$e(A/L\otimes_KA/J)$, where $L$ and $J$ denote left ideals of $A$;
moreover, an ``elimination lemma for quadric solvable polynomial
algebras" is obtained. By standard module theory, the obtained
results certainly apply to finitely generated $A$-modules.
{\parindent=0pt \v5

{\bf 3.1. Remark} Let $A$ be an {\it arbitrary} solvable polynomial
algebra $A$ with a  monomial ordering $\prec$ (which is
unnecessarily a graded monomial ordering). Note that
{\parindent=1.15truecm\par

\item{(1)} the calculation of Gelfand-Kirillov dimension and multiplicity
of a module $A/L$ deals only with the $\NZ$-filtration of $A/L$
induced by the standard $\NZ$-filtration of the $K$-vector space
$A$;

\item{(2)} with respect to {\it any} monomial ordering $\prec$ on $A$, $A$
has an $\prec$-filtration ${\cal F}A$ that turns $A$ into a
$\B$-filtered ring such that the associated graded algebra $G^{\cal
F}(A)$ is a homogeneous solvable polynomial algebra; and

\item{(3)} the tensor product $A\otimes_KB$ of $A$ with any solvable
polynomial algebras $B$ is a solvable polynomial algebra.
\par}

If one checks the text of ([Li1], CH.V, CH.VI) carefully, it is not
difficult to see that the same (key) result of ([Li], CH.V,
Proposition 7.2) holds for $A$, thereby the main results of ([Li1],
CH.V, CH.VI) all holds  true for $A$. \v5

{\bf 3.2. Remark} Furthermore, it was shown in ([Li1, LNM], CH.VIII,
Theorem 3.7, Theorem 4.1) that if $A$ is a quadric solvable
polynomial algebra with a graded monomial ordering $\prec_{gr}$,
then $A$ has global homological dimension gl.dim$A\le n$ and, by the
$K_0$-part of Quillen's theorem ([Qu], Theorem 7), the $K_0$-group
of $A$ is isomorphic to the additive group of integers $\mathbb{Z}$,
i.e., $K_0(A)\cong\mathbb{Z}$. Now, it follows from ([Li4], Ch.5,
Corollary 7.6), the $K_0$-part of Quillen's theorem ([Qu], Theorem
7), and our main result Theorem 2.5 that the same results hold true
in a more extended context, that is, we have\v5

{\bf 3.3. Proposition} Let $A=\KS /I$ be an algebra that satisfies
the equivalent conditions of Proposition 2.4. Then  gl.dim$A\le n$,
and  $K_0(A)\cong\mathbb{Z}$. Consequently, if $A=K[a_1,\ldots
,a_n]$ is  an arbitrary solvable polynomial algebra with a monomial
ordering $\prec$ (which is unnecessarily a graded monomial
ordering). Then gl.dim$A\le n$, and $K_0(A)\cong\mathbb{Z}$.}\v5

Before  giving a few observations on the algebras that satisfy the
equivalent conditions of Proposition 2.4, we write ${\cal A}_{\rm
gpbw}$ to denote the class of such algebras for convenience.
{\parindent=0pt\v5

{\bf Observation} (1)  If  $A=K\langle X_1,\ldots ,X_n\rangle /I$,
$B=K\langle Y_1,\ldots ,Y_m\rangle /J\in {\cal A}_{\rm gpbw}$, then,
by using an appropriate elimination monomial ordering, a similar
argument as in the proof (i) $\Rightarrow$ (ii) of Theorem 2.5
turns out that $A\otimes_KB\in {\cal A}_{\rm gpbw}$.
\par

(2) If  $A=\KS/I\in {\cal A}_{\rm gpbw}$, the Gr\"obner basis $\G=\{
g_{ji}~|~1\le i<j\le n\}$  of $I$ is obtained with respect to some
graded monomial ordering $\prec_{gr}$, and in $\G$ each
$g_{ji}=X_jX_i-\sum\mu_{k\ell}X_kX_{\ell}+\sum b_qX_q+c_{ji}$ with
$\mu_{k\ell},b_q, c_{ji}\in K$, and $X_kX_{\ell}\ne X_jX_i$, then,
it follows from ([Li4], Ch. 6, Theorem 3.1) that  $A$ is a
non-homogeneous Koszul algebra in the sense of [Pr] provided $b_q\in
K^*$ (if $b_q=0$, $c_{ji}=0$, then it is well-known that $A$ is a
homogeneous Koszul algebra). Consequently, every quadric solvable
polynomial algebra $A$ in the sense of [Li1] is either a homogeneous
Koszul algebra or a non-homogeneous Koszul algebra (see also [Li2],
Example 4.1).\par

(3) Let $A$ be as in (2) above. By ([Li4], Ch.7, Theorem 3.8), the
Rees algebra $\T A$ of $A$ defined with respect to the standard
$\NZ$-filtration of $A$ is a homogeneous Koszul algebra.}\v5

Thus, from both a computational and a structural viewpoint, our
discussion made so far has shown that it is worthwhile to pay more
attention to algebras in the class of algebras ${\cal A}_{\rm
gpbw}$. Especially for an algebra $A=\KS /I$ as described in Theorem
2.5(ii)(a), we have the following {\parindent=0pt\v5

{\bf Questions} (1) Is  $A$ a domain?\par

(2) Is $A$ a Noetherian ring?\par

(3) Is it always possible to define a monomial ordering $\prec$ on
the PBW $K$-basis $\B$ of $A$ via Proposition 2.6, such that $A$ is
turned into a solvable polynomial algebra?\v5

\centerline{References}
\parindent=1.3truecm

\item{[AL]} J. Apel and W. Lassner, An extension of Buchberger's
algorithm and calculations in enveloping fields of Lie algebras,
{\it J. Symbolic Comput}., 6(1988), 361--370.

\item{[Bu1]} B. Buchberger, {\it Ein Algorithmus zum Auffinden der
Basiselemente des Restklassenringes nach einem nulldimensionalen
polynomideal}, PhD thesis, University of Innsbruck, 1965.

\item{[Bu2]} B. Buchberger, Gr\"obner bases: ~An algorithmic method
in polynomial ideal theory. In: {\it Multidimensional Systems
Theory} (Bose, N.K., ed.), Reidel Dordrecht, 1985, 184--232.

\item{[BGV]} J. Bueso, J. G\'omez--Torrecillas, and A. Verschoren,
{\it Algorithmic methods in non-commutative algebra}: {\it
Applications to quantum groups}. Kluwer Academic Publishers, 2003.

\item{[EPS]} D. Eisenbud, I. Peeva and B. Sturmfels, Non-commutative
Gr\"obner bases for commutative algebras, {\it Proc. Amer. Math.
Soc.}, 126(1998), 687-691.

\item{[Gal]} A. Galligo, Some algorithmic questions on ideals of
differential operators, {\it Proc. EUROCAL'85}, LNCS 204, 1985,
413--421.

\item{[K-RW]} A. Kandri-Rody and V.~Weispfenning, Non-commutative
Gr\"obner bases in algebras of solvable type, {\it J. Symbolic
Comput.}, 9(1990), 1--26.

\item{[Li1]} H. Li, {\it Noncommutative Gr\"obner Bases and
Filtered-Graded Transfer}, LNM, 1795, Springer-Verlag, Berlin, 2002.

\item{[Li2]} H. Li, The general PBW ~property, {\it Alg.
Colloquium}, 14(4)(2007), 541--554.

\item{[Li3]} H. Li, Looking for Gr\"obner basis theory for (almost)
skew 2-nomial algebras, {\it J. Symbolic Computation}, 45(2010),
918--942.

\item{[Li4]} H. Li, {\it Gr\"obner Bases in Ring Theory}, World Scientific, 2011.

\item{[LW]} H. Li and  Y. Wu, Filtered-graded transfer of
Gr\"obner basis computation in solvable polynomial algebras, {\it
Comm. Alg.}, 28(1)(2000), 15--32.

\item{[Mor]} T. Mora, An introduction to commutative and noncommutative
Gr\"obner Bases, {\it Theoretic Computer Science}, 134(1994),
131--173.

\item{[Pr]} S. Priddy, Koszul resolutions, {\it Trans. Amer. Math.
Soc.}, 152(1970), 39--60.

\item{[Qu]} D. Quillen, Higher algebraic $K$-theory I, in {\it
Algebraic $K$-theory I}: Higher $K$-theory, ed., H. Bass, Lecture
Notes in Mathematics 341, Springer-Verlag, New York-Berlin, 1973,
85--147.

\end{document}